\documentclass[12pt]{amsart}

\begin{document}

\title{Genus Two Meromorphic Conformal Field Theory}
\author{Michael P. Tuite \\
Department of Mathematical Physics, National University of Ireland, Galway,
Ireland and The School of Theoretical Physics, Dublin Institute for Advanced
Studies, 10 Burlington Rd., Dublin 4, Ireland }
\date{Jan 10 2000}

\begin{abstract}
We construct the genus two (or two loop) partition function for meromorphic
bosonic conformal field theories. We use a sewing procedure involving two
genus one tori by exploiting an explicit relationship between the genus two
period matrix and pinching modular parameters. We obtain expressions for the
partition function for the chiral bosonic string, even rank lattice theories
and self-dual meromorphic conformal field theories including the Moonshine
Module. In particular, we find that for self-dual theories with central
charge 24, the genus two partition function multiplied by a universal
holomorphic function of the moduli is given by a meromorphic Siegel modular
form of weight $2$ where this universal function includes ghost
contributions. We also discuss a novel expansion for certain Siegel modular
forms.
\end{abstract}

\maketitle

\section{Introduction}

Since the earliest days of Conformal Field Theory (CFT) physicists have
studied multi-loop CFT on genus $g$ compact Riemann surfaces \cite{BPZ}, 
\cite{FS}, \cite{BK}, \cite{GSW}, \cite{P}. Most of this effort has been
devoted to the one loop case where a CFT is considered on a genus one torus
with a single modular parameter. Here the partition function can be
expressed geometrically by path integral methods in terms of formal operator
determinants or alternatively, using the quantum algebraic approach, in
terms of character trace formulas. For (chiral) bosonic Meromorphic CFTs
(MCFT) \cite{DGM} or Vertex Operator Algebras \cite{FLM}, \cite{Ka} (the
equivalent rigorous mathematical construction) the genus one partition
function is particularly simple and can be described by appropriate elliptic
modular forms \cite{DGM}, \cite{Z}. Examples include the chiral bosonic
string, even rank lattice theories and self-dual MCFTs such as the Moonshine
Module \cite{FLM}. For genus two or higher, most discussions have been of a
general nature and have ultimately expressed the partition function in terms
of various operator determinants e.g. \cite{BK}. In this paper, we
explicitly describe the partition function for bosonic Meromorphic CFTs on a
genus two Riemann surface. An expanded version of this work will appear
elsewhere.

We begin in Section 2 with a discussion of the moduli space for genus two
Riemann surfaces. Two separate modular parameterisations are of relevance:
(i) the genus two period matrix and (ii) `pinching'\ parameters used in a
standard formalism for sewing together two tori. Employing the genus two
period matrix, we describe the corresponding theory of Siegel modular forms
and discuss some basic Siegel forms and their Fourier expansions. With the
pinching parameters we can naturally describe the degenerations of a genus
two Riemann surface either into two tori or into a Riemann sphere.
Crucially, the two parameterisations can be related by an explicit formula.

In Section 3 we review some path integral results from classical string
theory. We discuss modular invariance and singular factorisation under
Riemann surface degenerations of the partition function for genus one and
genus two bosonic string theory. The importance of ghost contributions at
genus two is highlighted. In particular, one finds that the genus two
central charge $C=24$ bosonic string partition function when multiplied by a
certain universal holomorphic function of the moduli is given by the inverse
of a Siegel modular form of weight $10$ where the universal function
includes ghost contributions.

In Section 4 we discuss the genus two partition function for general bosonic
Meromorphic CFTs such as the chiral bosonic string, even self-dual lattice
theories and the Moonshine Module. We discuss an explicit formula for the
genus two partition function in terms of appropriate torus one point
functions. Employing the relationship between pinching parameters and the
period matrix allows us to express the partition function in terms of Siegel
modular forms. We find for self-dual theories with central charge 24 that
the genus two partition function multiplied by the universal holomorphic
function is a meromorphic Siegel form of weight $2$. We end with a number of
examples and a discussion of a novel expansion for certain Siegel forms in
terms of pinching parameters that explains observed patterns in their
Fourier coefficients.

\section{Genus Two Riemann Surfaces}

\subsection{Genus Two Moduli Space}

It is well known that the conformally inequivalent compact Riemann surfaces
of genus $g=1$ are described by one complex parameter whereas for $g\ge 2$,
they are described by $3g-3$ complex parameters forming the moduli space $%
\mathcal{M}_{g}$ \cite{FK}, \cite{IT}. $\mathcal{M}_{g}$ is analytic
everywhere except at a discrete set of points and has a universal covering
space, known as the Teichm\"{u}ller space $\mathcal{T}_{g}$, which is a
smooth complex analytic manifold. Then $\mathcal{M}_{g}$=$\mathcal{T}%
_{g}/\Gamma $ where $\Gamma $ is the mapping class group, a discrete group
acting on $\mathcal{T}_{g}$. Genus one Riemann surfaces are parameterised by
the standard modular parameter $\tau \in \mathbb{H}$, the upper half complex
plane with $Im(\tau )>0$ and $\Gamma =PSL(2,\mathbb{Z}).$ For a genus two
Riemann surface, the Teichm\"{u}ller space $\mathcal{T}_{2}$ is of complex
dimension three and can be described in various ways. Two separate
descriptions are of particular interest to us, namely, (i) the genus two
period matrix and (ii) so-called `pinching'\ parameters used in a standard
formalism for sewing together Riemann surfaces. In this section we discuss
these two parameterisations and how they are explicitly related. We also
consider the degeneration properties of genus two Riemann surfaces and
lastly discuss relevant aspects of the theory of Siegel modular forms.

\subsection{Holomorphic 1 Forms and Some Genus One Theory}

Consider a compact Riemann surface $\mathcal{S}$ of genus $g$ and consider
the canonical homology basis $a_{1,}\ldots a_{g},b_{1},\ldots b_{g}$ which
is illustrated in Fig. 1 for a genus two Riemann surface. In general, there
exists a set of $g$ holomorphic 1-forms $\nu _{i}^{(g)}$, $i=1,\ldots g$
normalised by \cite{FK}, \cite{Sp} 
\begin{equation}
\int_{a_{i}}\nu _{j}^{(g)}=\delta _{ij}{.}  \label{norm}
\end{equation}
%
%
%
%
%
%
\vskip 0.5 truecm
\begin{center}
\begin{picture}(200,100)
\put(0,50){\qbezier(10,18)(50,35)(90,18)} \put(0,50){\qbezier(10,18)(-20,0)(10,-18)} \put(0,50){\qbezier(10,-18)(50,-35)(90,-18)} 
 \put(0,50){\qbezier(90,18)(100,12)(110,18)} \put(0,50){\qbezier(90,-18)(100,-12)(110,-18)} 
\put(0,50){\qbezier(25,0)(50,17)(75,0)} \put(0,50){\qbezier(20,2)(50,-17)(80,2)} \put(100,50){\qbezier(10,18)(50,35)(90,18)} \put(100,50){\qbezier(10,-18)(50,-35)(90,-18)} \put(100,50){\qbezier(90,-18)(120,0)(90,18)} 
\put(100,50){\qbezier(25,0)(50,17)(75,0)} \put(100,50){\qbezier(20,2)(50,-17)(80,2)} 
 \put(0,50){\qbezier(50,-26)(44,-17)(50,-8)} \put(0,50){\qbezier[15](50,-26)(56,-17)(50,-8)} \put(47,34){\vector(0,1){0}} \put(41,38){\makebox(0,0){$a_1$}}

\put(100,50){\qbezier(50,-26)(44,-17)(50,-8)} \put(100,50){\qbezier[15](50,-26)(56,-17)(50,-8)} \put(147,34){\vector(0,1){0}} \put(141,38){\makebox(0,0){$a_2$}}

\put(0,50){\qbezier(17,9)(50,26)(82,9)} \put(0,50){\qbezier(17,9)(0,0)(17,-9)} \put(0,50){\qbezier(17,-9)(50,-26)(82,-9)} \put(0,50){\qbezier(82,-9)(100,0)(82,9)} \put(9,48){\vector(0,-1){0}} \put(3,50){\makebox(0,0){$b_1$}}

\put(100,50){\qbezier(17,9)(50,26)(82,9)} \put(100,50){\qbezier(17,9)(0,0)(17,-9)} \put(100,50){\qbezier(17,-9)(50,-26)(82,-9)} \put(100,50){\qbezier(82,-9)(100,0)(82,9)} \put(109,48){\vector(0,-1){0}} \put(103,50){\makebox(0,0){$b_2$}}

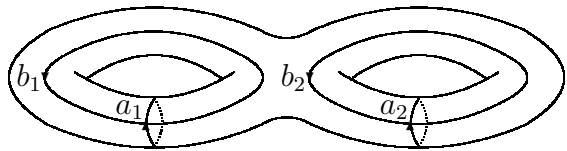
\put(100,0){\makebox(0,0){\footnotesize{Fig. 1. The Genus Two Homology Basis}}}
\end{picture}

\end{center}
\vskip 0.5 truecm
%
These 1-forms can be neatly encapsulated in a single singular bilinear 2
form $\omega ^{(g)}$, known as the normalised differential of the second
kind, which is defined by the following properties \cite{Sp}, \cite{Y}: 
\begin{equation}
\omega ^{(g)}(x,y)=(\frac{1}{(x-y)^{2}}+\text{regular terms})dxdy
\label{omega}
\end{equation}
for $x,y\in \mathcal{S}$ with normalisation $\int_{a_{i}}\omega ^{(g)}(\cdot
,y)=0${\ for }$i=1,\ldots g$. Then, using the Riemann bilinear relations,
one finds that

\begin{equation}
\nu _{i}^{(g)}=\frac{1}{2\pi i}\int_{b_{i}}\omega ^{(g)}(\cdot ,y)
\label{nui}
\end{equation}
with $\nu _{i}^{(g)}$normalised as in (\ref{norm}).

As an example, consider the standard genus one Riemann torus with periods $1$
and $\tau $ along the $a$ and $b$ cycles. Then the normalised holomorphic
1-form is $\nu ^{(1)}=dz\,$ and the normalised differential of the second
kind is 
\begin{equation}
\omega ^{(1)}(x,y)=(\wp (x-y)+\hat{E}_{2})dxdy,  \label{2form}
\end{equation}
where 
\begin{equation}
\wp (z)=\frac{1}{z^{2}}+\sum_{k=2}^{\infty }\hat{E}_{2k}z^{2k-2},
\end{equation}
is the Weierstrass function \cite{FK}, \cite{Sp} and 
\begin{eqnarray}
\hat{E}_{2k}(q) &=&-\frac{B_{2k}}{(2k)!}E_{2k}(q)=-\frac{B_{2k}}{(2k)!}+%
\frac{4k}{(2k)!}q+O(q^{2}),  \label{ehat} \\
E_{2k}(q) &=&1-\frac{4k}{B_{2k}}\sum_{n=1}^{\infty }\sigma _{2k-1}(n)q^{n}.
\label{elleisen}
\end{eqnarray}
Here $E_{2k}(q)$ is the standard elliptic Eisenstein series, $q=\exp (2\pi
i\tau )$, $\sigma _{2k-1}(n)=\sum_{d\mid n}d^{2k-1}$ and $B_{2k}$ is the $2k$%
th Bernoulli number defined by 
\begin{equation}
\frac{t}{e^{t}-1}=\sum_{k=0}^{\infty }B_{2k}\frac{t^{2k}}{(2k)!}.
\end{equation}
$\hat{E}_{2k}(q)$ is a $PSL(2,\mathbb{Z)}$ modular form of weight $2k$ for $%
k>1$. $\hat{E}_{2}(q)$ transforms anomalously under $PSL(2,\mathbb{Z)}$
where in particular \cite{Se}, \cite{DLM} 
\begin{equation}
\hat{E}_{2}\rightarrow \tau ^{2}(\hat{E}_{2}-\frac{1}{2\pi i\tau }).
\label{SEhat}
\end{equation}
under $S:\tau \rightarrow -1/\tau $. $\hat{E}_{2}$ plays a fundamental role
throughout our analysis and also appears in the modular 'covariant'
derivative for elliptic modular forms defined as follows. Let $f_{k}(q)$ be
a $PSL(2,\mathbb{Z)}$ modular form of weight $k$. Then 
\begin{equation}
D_{q}f_{k}\equiv \left( q\frac{d}{dq}+k\hat{E}_{2}(q)\right) f_{k}
\label{Df}
\end{equation}
is a modular form of weight $k+2$. This covariant derivative will also
appear later on.

\subsection{The Period Matrix and Symplectic Modular Group}

The genus $g$ period matrix $\Omega $ is defined by \cite{FK}, \cite{Sp} 
\begin{equation}
\Omega _{ij}=\int_{b_{i}}\nu _{j}^{(g)},\quad i,j=1\ldots g{.}  \label{Omega}
\end{equation}
It is well known that $\Omega $ is a symmetric complex matrix with positive
imaginary part so that $\Omega \in \mathbb{H}_{g}$, the $g$ dimensional
Siegel complex upper half plane. There is a natural symplectic intersection
form $\Xi $ acting on the homology basis with 
\begin{equation}
\Xi (a_{i},a_{j})=\Xi (b_{i},b_{j})=0,\quad \Xi (a_{i},b_{j})=\delta
_{ij},\quad i,j=1\ldots g
\end{equation}
which is preserved by the symplectic group $Sp(2g,\mathbb{Z})$ acting on the
homology basis with $\gamma ^{T}\Xi \gamma =\Xi $ for each $\gamma \in Sp(2g,%
\mathbb{Z})$ where 
\begin{equation}
{\gamma =}\left( 
\begin{array}{ll}
A & B \\ 
C & D
\end{array}
\right) {,\quad \gamma ^{-1}=}\left( {\ 
\begin{array}{ll}
{{D^{T}}} & {{-B^{T}}} \\ 
{{-C^{T}}} & {{A^{T}}}
\end{array}
}\right)  \label{sp4Z}
\end{equation}
and where $A,B,C,D$ are integer valued $g\times g$ matrices \cite{FK}, \cite
{Fr}, \cite{Kl}. The action of $\gamma $ on the homology basis induces a
modular transformation on $\Omega $ 
\begin{equation}
{\gamma [\Omega ]=(A\Omega +B)(C\Omega +D)^{-1},}  \label{modtrans}
\end{equation}
where $\gamma =\pm I_{2g}$ acts as the identity.

For a genus two Riemann surface, the elements of the period matrix form a
set of 3 independent parameters so that the genus two Teichm\"{u}ller space
is $\mathcal{T}_{2}=\mathbb{H}_{2}$ \cite{FK}. The modular group acting on $%
\mathcal{T}_{2}$ is then given by $PSp(4,\mathbb{Z})\equiv Sp(4,\mathbb{Z}%
)/\langle \pm I_{4}\rangle $ so that the genus two moduli space is $\mathcal{%
M}_{2}=\mathbb{H}_{2}/PSp(4,\mathbb{Z})$. It is useful to give an explicit
basis of generators for $\Gamma $ as follows: 
\begin{eqnarray}  \label{genset}
\\
S_{1} &=&\left[ 
\begin{array}{llll}
0 & 0 & 1 & 0 \\ 
0 & 1 & 0 & 0 \\ 
-1 & 0 & 0 & 0 \\ 
0 & 0 & 0 & 1
\end{array}
\right], S_{2}=\left[ 
\begin{array}{llll}
1 & 0 & 0 & 0 \\ 
0 & 0 & 0 & 1 \\ 
0 & 0 & 1 & 0 \\ 
0 & -1 & 0 & 0
\end{array}
\right],  \nonumber \\
T_{1}&=&\left[ 
\begin{array}{llll}
1 & 0 & 1 & 0 \\ 
0 & 1 & 0 & 0 \\ 
0 & 0 & 1 & 0 \\ 
0 & 0 & 0 & 1
\end{array}
\right], T_{2}=\left[ 
\begin{array}{llll}
1 & 0 & 0 & 0 \\ 
0 & 1 & 0 & 1 \\ 
0 & 0 & 1 & 0 \\ 
0 & 0 & 0 & 1
\end{array}
\right], U=\left[ 
\begin{array}{llll}
1 & 0 & 0 & 1 \\ 
0 & 1 & 1 & 0 \\ 
0 & 0 & 1 & 0 \\ 
0 & 0 & 0 & 1
\end{array}
\right].  \nonumber
\end{eqnarray}
The generators $T_{1}$, $T_{2}$ and $U$ define the translations $\Omega
_{11}\rightarrow \Omega _{11}+1$, $\Omega _{22}\rightarrow \Omega _{22}+1$
and $\Omega _{12}\rightarrow \Omega _{12}+1$ respectively and $V\equiv 
\mathrm{diag}(-1,1,-1,1)$ ($V=S_{1}^{2}=-S_{2}^{2}$) describes the
reflection $\Omega _{12}\rightarrow -\Omega _{12}$. The subgroups $\langle
S_{1},T_{1}\rangle $ and $\langle S_{2},T_{2}\rangle $ are isomorphic to $%
PSL(2,\mathbb{Z})$ and correspond to the modular groups for the left and
right hand tori (see below). The period matrix parameterisation of moduli
space is particularly suitable for constructing modular forms on $\mathcal{M}%
_{2}$. However, an alternative parameterisation is more suitable for
describing degenerations of a genus two Riemann surface as we now explain.

\subsection{Pinching Parameters and Riemann Surface Degenerations}

An alternative way of parameterising $\mathcal{M}_{2}$ is provided by a
standard sewing procedure for joining together two Riemann surfaces or for
attaching a handle to a Riemann surface \cite{Kn}, \cite{So}, \cite{IT}. In
general, consider two compact Riemann surfaces $\mathcal{S}_{1}$ and $%
\mathcal{S}_{2}$ of genus $g_{1}$ and $g_{2}$. These can be sewn together to
form a Riemann surface of genus $g_{1}+g_{2}$ as follows. Choose local
coordinates $z_{1}$ on $\mathcal{S}_{1}$ and $z_{2}$ on $\mathcal{S}_{2}$
and excise the two disks $\left| z_{1}\right| $ $<|\epsilon |$ and $\left|
z_{2}\right| <|\epsilon |$ where $\epsilon $ is a complex parameter with $%
|\epsilon |<1$. The two surfaces are sewn together by identifying the
annular regions $|\epsilon |\le \left| z_{1}\right| \le 1$ and $|\epsilon
|\le \left| z_{2}\right| \le 1$ via the relation 
\begin{equation}
z_{1}z_{2}={\epsilon }  \label{pinch}
\end{equation}
\vskip 0.5 truecm
\begin{center}
\begin{picture}(250,100)
\put(0,50){\qbezier(10,18)(50,35)(90,18)} \put(0,50){\qbezier(10,18)(-20,0)(10,-18)} \put(0,50){\qbezier(10,-18)(50,-35)(90,-18)} \put(0,50){\qbezier(25,0)(50,17)(75,0)} \put(0,50){\qbezier(20,2)(50,-17)(80,2)} 
\put(125,50){\qbezier(10,18)(50,35)(90,18)} \put(125,50){\qbezier(10,-18)(50,-35)(90,-18)} \put(125,50){\qbezier(90,-18)(120,0)(90,18)} 
\put(125,50){\qbezier(25,0)(50,17)(75,0)} \put(125,50){\qbezier(20,2)(50,-17)(80,2)} 
  \put(0,50){\qbezier(90,18)(84,18)(84,0)} \put(0,50){\qbezier(84,0)(84,-18)(90,-18)} \put(0,50){\qbezier(90,18)(96,18)(96,0)} \put(0,50){\qbezier(96,0)(96,-18)(90,-18)} 
 \put(0,50){\qbezier(90,12)(87,12)(87,0)} \put(0,50){\qbezier(87,0)(87,-12)(90,-12)} \put(0,50){\qbezier(90,12)(93,12)(93,0)} \put(0,50){\qbezier(93,0)(93,-12)(90,-12)} 
 \put(45,50){\qbezier(90,18)(84,18)(84,0)} \put(45,50){\qbezier(84,0)(84,-18)(90,-18)} \put(45,50){\qbezier(90,18)(96,18)(96,0)} \put(45,50){\qbezier(96,0)(96,-18)(90,-18)} 
 \put(45,50){\qbezier(90,12)(87,12)(87,0)} \put(45,50){\qbezier(87,0)(87,-12)(90,-12)} \put(45,50){\qbezier(90,12)(93,12)(93,0)} \put(45,50){\qbezier(93,0)(93,-12)(90,-12)} 
 
\put(0,50){\qbezier(90,18)(112,15)(135,12)} \put(0,50){\qbezier(90,-18)(112,-15)(135,-12)} 
\put(0,50){\qbezier(90,12)[10](112,15)(135,18)} \put(0,50){\qbezier(90,-12)[10](112,-15)(135,-18)} 
 \put(90,50){\vector(-1,-2){0}} \put(0,50){\qbezier(90,0)(100,15)(90,30)} 
\put(90,90){\makebox(0,0){$z_1=0$}}

\put(135,50){\vector(1,-2){0}} \put(45,50){\qbezier(90,0)(80,15)(90,30)} 
\put(135,90){\makebox(0,0){$z_2=0$}}

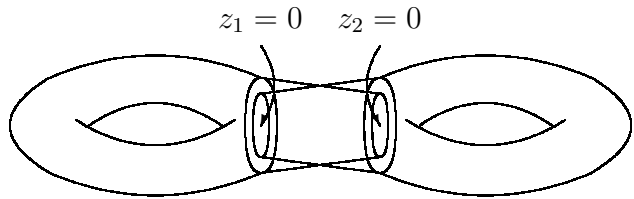
\put(100,0){\makebox(0,0){\footnotesize{Fig. 2. Sewing Two Tori}}}
\end{picture}
\end{center}
\vskip 0.5 truecm
%
Thus every genus two Riemann surface can be constructed by sewing together
two tori as is depicted in Fig. 2 where the identified annular regions can
be thought of as a connecting cylinder. \noindent Note that as $\epsilon
\rightarrow 0$, the cylinder is pinched down and the Riemann surface
degenerates into two tori with standard modular parameters $\tau _{1}$ and $%
\tau _{2}$. Furthermore, the parameters $q_{1}=\exp (2\pi i\tau _{1})$ and $%
q_{2}=\exp (2\pi i\tau _{2})$ can themselves be interpreted as the sewing
parameters associated with gluing two disks of a Riemann sphere together to
form a torus where the torus degenerates to a Riemann sphere as $%
q\rightarrow 0$. The parameters $q_{1},q_{2},\epsilon $ therefore form an
alternative set of modular parameters for genus two Riemann surfaces.

These parameters are particularly suitable for describing all degenerations
of the genus two Riemann surface. In particular, the degenerations into one
torus with modular parameter $q_{1}$ (or $q_{2}$) described by $%
q_{2}\rightarrow 0$ (or $q_{1}\rightarrow 0$) or into two tori described by $%
\epsilon \rightarrow 0$ are naturally described. \noindent The degeneration
into two tori where $\epsilon \rightarrow 0$ is of particular importance in
our subsequent discussion.

\subsection{Relating the Two Parameterisations}

The genus two period matrix ${\Omega }$ is a holomorphic function of the
pinching parameters $q_{1},q_{2},\epsilon $. A general method for
calculating ${\Omega }$ for any two sewn Riemann surfaces has been described
by Yamada \cite{Y}. Using this formalism, the genus two normalised
differential of the second kind $\omega ^{(2)}\,$of (\ref{omega}) can be
expressed as a power series in $\epsilon $ with coefficients calculated from
the genus one normalised differential $\omega ^{(1)}$ of (\ref{2form}) for
the left and right torus. Integrating $\omega ^{(2)}$ over the appropriate
cycles as in (\ref{nui}) and (\ref{Omega}) eventually results in the
following exact expressions for the period matrix elements: 
\begin{eqnarray}
{\Omega }_{11} &=&\tau _{1}+\frac{\epsilon }{2\pi i}\left(
A(q_{2})(1-A(q_{1})A(q_{2}))^{-1}\right) _{11}  \label{om11} \\
&=&\tau _{1}+\frac{\epsilon ^{2}}{2\pi i}\hat{E}_{2}(q_{2})+O(\epsilon ^{4}),
\nonumber \\
{\Omega }_{12} &=&-\frac{\epsilon }{2\pi i}\left( 1-A(q_{1})A(q_{2})\right)
_{11}^{-1}  \label{om12} \\
&=&-\frac{\epsilon }{2\pi i}\left( 1+\hat{E}_{2}(q_{1})\hat{E}%
_{2}(q_{2})\epsilon ^{2}\right) +O(\epsilon ^{5}),  \nonumber
\end{eqnarray}
with ${\Omega }_{22}={\Omega }_{11}(q_{1}\leftrightarrow q_{2})$. The
subscript $11$ appearing on the right hand side of (\ref{om11}) and (\ref
{om12}) denotes the first component of matrices constructed from the
infinite dimensional matrix $A$ with components 
\begin{equation}
A_{mn}(q)=\epsilon ^{m+n-1}\binom{2m+2n-3}{2m-1}\hat{E}_{2m+2n-2}(q)
\label{A}
\end{equation}
using (\ref{ehat}). It is not difficult to check that ${\Omega }_{ij}$
transforms correctly under the action of the generators $\{S_{1},T_{1}\}$
given in (\ref{genset}) of the left $PSL(2,\mathbb{Z)}$ modular subgroup of $%
{Sp}(4,\mathbb{Z})$ where for 
\begin{equation}
S_{1}:\tau _{1}\rightarrow -1/\tau _{1},\ \epsilon \rightarrow -\epsilon
/\tau _{1},  \label{S1}
\end{equation}
using the anomalous modular transformation property for $\hat{E}_{2}(q)\,$in
(\ref{SEhat}), one finds that 
\begin{eqnarray}
{\Omega }_{11} &\rightarrow &-1/{\Omega }_{11},\,  \nonumber \\
{\Omega }_{12} &\rightarrow &-{\Omega }_{12}/{\Omega }_{11},\,  \nonumber \\
{\Omega }_{22} &\rightarrow &{\Omega }_{22}-{\Omega }_{12}^{2}/{\Omega }%
_{11},  \label{SOmega}
\end{eqnarray}
as expected from (\ref{modtrans}). A similar result hold for $S_{2}$. Notice
also that the transformation $\epsilon \rightarrow -\epsilon $ is equivalent
to the reflection symmetry $V:{\Omega }_{12}\rightarrow -{\Omega }_{12}.$

\subsection{$P{Sp}(4,\mathbb{Z})$ Siegel Modular Forms}

A $PSp(4,\mathbb{Z})\,$Siegel modular form of weight $k$ is a holomorphic
function on $\mathbb{H}_{2}$ satisfying \cite{Ig}, \cite{Fr} 
\begin{equation}
f(\gamma [\Omega ]){=\mathrm{det}}(C\Omega +D)^{k}f(\Omega ){,}
\label{modform}
\end{equation}
where $\gamma \in PSp(4,\mathbb{Z})$, as given in (\ref{sp4Z}) for $g=2$. A
cusp form is one which vanishes under any Riemann surface degeneration.
These definitions can be extended to include forms with multiplier systems.
We will also refer to any meromorphic function which transforms as in (\ref
{modform}) as a meromorphic modular form.

Every genus two Siegel modular form is invariant under the translations $%
T_{1},T_{2}$ and $U$ of (\ref{genset}) since $\mathrm{{det}(C\Omega +D)=1}$
in each case. Therefore every modular form can be expanded as a Fourier
series in the following three variables: 
\begin{equation}
q=e^{2\pi i\Omega _{11}},\quad r=e^{2\pi i\Omega _{12}},\quad s=e^{2\pi
i\Omega _{22}}.  \label{qrs}
\end{equation}
It can be shown for every $PSp(4,\mathbb{Z})\,$ Siegel form that the Fourier
coefficient of $q^{a}r^{b}s^{c}$ is non-vanishing provided $a,c\geq 0$ and $%
b^{2}\leq 4ac$ \cite{Ig}. Furthermore, it is useful for us to introduce 
\begin{equation}
u=r+r^{-1}-2,  \label{u}
\end{equation}
which is invariant under the reflection $V$.

The Fourier parameters $q,s,u$ can be related to the pinching parameters via
(\ref{om11}) and (\ref{om12}) to obtain: 
\begin{eqnarray}
q &=&q_{1}(1+\epsilon ^{2}\hat{E}_{2}(q_{2}))+O(\epsilon ^{4}),
\label{qsueps} \\
s&=&q_{2}(1+\epsilon ^{2}\hat{E}_{2}(q_{1}))+O(\epsilon ^{4}),  \nonumber \\
u &=&\epsilon ^{2}(1+\epsilon ^{2}(\frac{1}{12}+2\hat{E}_{2}(q_{1})\hat{E}%
_{2}(q_{2}))+O(\epsilon ^{6}).  \nonumber
\end{eqnarray}
Note that the degeneration $q_{1}\rightarrow 0$ ($q_{2}\rightarrow 0$) of
the genus two Riemann surface to the right (left) torus is given by $%
q\rightarrow 0$ ($s\rightarrow 0$) whereas the degeneration $\epsilon
\rightarrow 0$ to two tori is given by $u\rightarrow 0$ with these
parameters.

Three standard sets of Siegel modular forms will be useful to us:
Half-integral Characteristic Theta Series and Lattice Theta Series and
Eisenstein Series.

\textbf{Half-integral Characteristic Theta Series. }These are defined by 
\begin{eqnarray}
{\Theta }\left[ 
\begin{array}{l}
\mathbf{a} \\ 
\mathbf{b}
\end{array}
\right] {(\Omega )} &=&{}{{\sum_{\mathbf{n}\in \mathbb{Z\times Z}}\exp (}}%
i\pi (\mathbf{n}+\mathbf{a}).\Omega .(\mathbf{n}+\mathbf{a})+2\pi i(\mathbf{n%
}+\mathbf{a}).\mathbf{b}{{\mathbf{),}}}  \nonumber \\
{{\mathbf{a,b}}} &\in &{}{{\mathbf{\{}}}\left[ 
\begin{array}{l}
0 \\ 
0
\end{array}
\right] ,\left[ 
\begin{array}{l}
0 \\ 
\frac{1}{2}
\end{array}
\right] ,\left[ 
\begin{array}{l}
\frac{1}{2} \\ 
0
\end{array}
\right] ,\left[ 
\begin{array}{l}
\frac{1}{2} \\ 
\frac{1}{2}
\end{array}
\right] \}.  \label{theta}
\end{eqnarray}
There are 10 independent non-zero theta series which transform amongst
themselves with weight $\frac{1}{2}$ and multiplier system \cite{Mu}, \cite
{Fr}, \cite{Kl}. Under the degeneration $q_{1}\rightarrow 0$, ${\Theta }%
\left[ 
\begin{array}{l}
\mathbf{a} \\ 
\mathbf{b}
\end{array}
\right] {(\Omega )}\rightarrow \theta \left[ 
\begin{array}{l}
a_{2} \\ 
b_{2}
\end{array}
\right] (\tau _{2})$ where $\theta \left[ 
\begin{array}{l}
a \\ 
b
\end{array}
\right]$ is one the three standard elliptic theta series with $a,b\in \{0,%
\frac{1}{2}\}$ (and similarly for $q_{2}\rightarrow 0$). For $\epsilon
\rightarrow 0,$ ${\Theta }\left[ 
\begin{array}{l}
\mathbf{a} \\ 
\mathbf{b}
\end{array}
\right] {(\Omega )}\rightarrow \theta \left[ 
\begin{array}{l}
a_{1} \\ 
b_{1}
\end{array}
\right] (\tau _{1})\theta \left[ 
\begin{array}{l}
a_{2} \\ 
b_{2}
\end{array}
\right] (\tau _{2})$.

\textbf{Lattice Theta Series.} Let $\Lambda $ be an integral self dual
lattice with Euclidean metric $\langle ,\rangle $ of rank $n=0\text{ mod }8$
. Then 
\begin{equation}
\Theta _{\Lambda }{(\Omega )}={{\sum_{\alpha ,\beta \in \Lambda }}}%
q^{\langle \alpha ,\alpha \rangle /2}s^{\langle \beta ,\beta \rangle
/2}r^{\langle \alpha ,\beta \rangle }\mathbf{,}  \label{ThetaL}
\end{equation}
is a Siegel modular form of weight $n/2.$ Under the degeneration $%
q_{1}\rightarrow 0$, $\Theta _{\Lambda }{(\Omega )}\rightarrow \theta
_{\Lambda }(q_{2})$ and (similarly for $q_{2}\rightarrow 0$) whereas under ${%
\epsilon }\rightarrow 0$, $\Theta _{\Lambda }{(\Omega )}\rightarrow \theta
_{\Lambda }{(}q_{1}{)}\theta _{\Lambda }(q_{2})$ where $\theta _{\Lambda
}(q)=\sum_{\alpha \in \Lambda }q^{\langle \alpha ,\alpha \rangle /2}$ is the
standard elliptic lattice theta function.

\textbf{Eisenstein Series. }The genus two Eisenstein series for $P{Sp}(4,%
\mathbb{Z})$ is defined by 
\begin{equation}
\psi _{2k}(\Omega )=\sum_{\{C,D\}}{(\mathrm{det}(C\Omega +D))}^{-2k}{,}
\label{eisen}
\end{equation}
where the sum is performed over all inequivalent pairs of matrices $C,D$
forming the lower two rows of the elements of $PSp(4,\mathbb{Z})$ \cite{Fr}, 
\cite{Kl}. For any integer $k>1$, $\psi _{2k}(\Omega )$ is a Siegel modular
form of weight $2k$. For $q_{1}\rightarrow 0$ we have $\psi _{2k}(\Omega
)\rightarrow E_{2k}(q_{2})$ with $E_{2k}$ as in (\ref{ehat}) (and similarly
for $q_{2}\rightarrow 0$). On the other hand, for $\epsilon \rightarrow 0$, $%
\psi _{2k}(\Omega )\rightarrow \phi _{2k}(q_{1},q_{2})$, an elliptic modular 
$2k$ form for both $q_{1}$ and $q_{2}$ and symmetric in its two arguments.
Elliptic forms are unique for weight $2k\leq 10$ so that $\phi
_{2k}(q_{1},q_{2})=E_{2k}(q_{1})E_{2k}(q_{2})$ in these cases.

The first few coefficients of the Fourier expansion for $\psi _{4}$, $\psi
_{6}$ in terms of $q,s,u$ are: 
\begin{eqnarray}
\psi _{4} &=& E_4(q)E_4(s) +14400qsu+240qsu^{2}+O(q^{2},s^{2}),  \label{psi4}
\\
\psi _{6} &=& E_6(q)E_6(s) +42336qsu-504qsu^{2}+O(q^{2},s^{2}),  \label{psi6}
\end{eqnarray}
where $E_{4}(q)=1+240q+O(q^{2})$ and $E_{6}(q)=1-504q+O(q^{2})$. We note
that the $qsu^{2}$ (or $qsr^{2},qsr^{-2}$ ) term has been mistakenly omitted
from the Fourier expansions given by Igusa in \cite{Ig} and quoted in \cite
{Kl}. We obtained these Fourier coefficients from expressions for these
forms in terms of the half integral theta series. We also observe that $%
14400=(240)^{2}/4$ and $42336=(-504)^{2}/6$.

The ring of all Siegel modular forms for $PSp(4,\mathbb{Z})$ has been
determined by Igusa to be generated by a set of four independent forms of
weight 4, 6, 10 and 12 \cite{Ig}. Using the Eisenstein and Theta series
these can be explicitly chosen to be $\psi _{4}$, $\psi _{6}$, $\Delta _{10}$
and $F_{12}$ with

\begin{eqnarray}
\Delta _{10}(\Omega ) &=&\frac{1}{2^{12}}\prod_{\mathbf{a},\mathbf{b}}{%
\Theta }\left[ 
\begin{array}{l}
\mathbf{a} \\ 
\mathbf{b}
\end{array}
\right] ^{2},  \label{D10} \\
F_{12}(\Omega ) &=&{\frac{1}{4}}\sum_{\mathbf{a},\mathbf{b}}{\Theta }\left[ 
\begin{array}{l}
\mathbf{a} \\ 
\mathbf{b}
\end{array}
\right] ^{24},  \label{F12}
\end{eqnarray}
where the sum and product are taken over the ten independent theta series.
The first few coefficients of the Fourier expansion are: 
\begin{eqnarray}
\Delta _{10} &=&qsu((1-24q)(1-24s)  \label{D10fourier} \\
&&+2u(-q-s+72qs)-16qsu^{2})+O(q^{3},s^{3}),  \nonumber \\
F_{12} &=& f_{12}(q)f_{12}(s) +101568qsu+1104qsu^{2}+O(q^{2},s^{2}),
\label{F12fourier}
\end{eqnarray}
with $f_{12}(q)=\sum {\theta }^{24}(q)/2=1+1104q+O(q^{2})$ where the sum is
taken over the three elliptic theta functions.

$\Delta _{10}$ is a Siegel form of weight 10 which plays a central role in
genus two Meromorphic CFT. $\Delta _{10}$ is the unique cusp form of weight
10 and vanishes for all degenerations of the genus two surface. Under the
degeneration, $\epsilon \rightarrow 0,$ we have the following important
factorisation property using (\ref{qsueps}) 
\begin{equation}
\Delta _{10}(\Omega )=\epsilon ^{2}\Delta (q_{1})\Delta (q_{2})+O(\epsilon
^{4}),  \label{D10deg}
\end{equation}
where $\Delta (q)=\eta ^{24}(q)$ is the elliptic cusp form of weight 12 and $%
\eta (q)=q^{1/24}\prod_{n=1}^{\infty }(1-q^{n})$ is the Dedekind eta.  The
existence of the cusp form $\Delta _{10}$ at genus two is related that for $%
\Delta (q)$ at genus one which vanishes as the torus degenerates with $%
q\rightarrow 0$ as is further discussed below \cite{Kn}, \cite{Mo}, \cite
{KMO}.

Under the degeneration, $\epsilon \rightarrow 0,$ the 12 form $F_{12}$
factorises with $F_{12}$ $\rightarrow f_{12}(q_{1})f_{12}(q_{2})$. We
observe that $101568=(1104)^{2}/12$ which leads us to the following
observation: Let $F_{k}(\Omega )$ be a Siegel form of weight $k=4,6,8$ and $%
12$ which factorises under $\epsilon \rightarrow 0$ with $F_{k}$ $%
\rightarrow f_{k}(q_{1})f_{k}(q_{2})$ for some elliptic modular $k$ form $%
f_{k}(q)=1+aq+O(q^{2})$ then 
\begin{equation}
F_{k}(\Omega )=(1+aq)(1+as)+\frac{a^{2}}{k}qsu+aqsu^{2}+O(q^{2},s^{2}).
\label{FkFourier}
\end{equation}
For a self-dual lattice $\Lambda \,$of rank $2k$ $=0\text{ mod }8$ we have $%
F_{k}(\Omega )=\Theta _{\Lambda }(\Omega )$. Then $a$ is the multiplicity of
lattice roots and is also the coefficient of $qsu^{2}$ being the
multiplicity of pairs of roots $\alpha ,\beta $ such that $\langle \alpha
,\beta \rangle =2$ i.e. $\alpha =\beta $. We will return to a CFT
explanation for the $\frac{a^{2}}{k}$coefficient of $qsu$ later on for
weights $k=4,8,12$.

Finally, there are three independent 12 forms which can be conveniently
chosen to be $F_{12}$, $\psi _{4}^{3}$ and $\psi _{6}^{2}$. It is a very
important observation for our purposes to note that there are no forms of
weight 12 which vanish under the degeneration $\epsilon \rightarrow 0$ since
there are three independent symmetric elliptic 12 forms in $q_{1}$ and $%
q_{2} $ i.e. $E_{12}(q_{1})E_{12}(q_{2})$, $E_{12}(q_{1})\Delta
(q_{2})+\Delta (q_{1})E_{12}(q_{2})$ and $\Delta (q_{1})\Delta (q_{2}).$
This implies that every $PSp(4,\mathbb{Z})$ Siegel 12 form is uniquely
determined by its degeneration under $\epsilon \rightarrow 0.$ Thus, for
every Niemeier lattice $\Lambda ,$ $\Theta _{\Lambda }(\Omega )$ is the
unique 12 form for which $\Theta _{\Lambda }(\Omega )\rightarrow \theta
_{\Lambda }(q_{1})\theta _{\Lambda }(q_{2})$ under $\epsilon \rightarrow 0.$

\section{The Genus Two Polyakov String}

We now briefly review some relevant classical results for the 26 dimensional
Polyakov bosonic string. Here the multi-loop probability amplitude has been
known for quite some time and is found by path integral methods in terms of
various determinants \cite{GSW}, \cite{BK}, \cite{Kn}. The purpose of this
section is to discuss the origin of modular invariance, factorisation and
the singularity structure of the amplitude in this geometrical approach. We
will identify the contributions to the one and two loop probability
amplitude which arise from the partition function for the underlying bosonic
central charge $C=24$ Meromorphic Conformal Field Theory (MCFT) for genus
one and genus two. Despite the fact that the methods described in this
section may not yet have a rigorous mathematical formulation, the insights
gained are invaluable. In particular, we wish to highlight the important
role played by ghost contributions at genus two in order to achieve full
modular invariance.

In 26 dimensional Polyakov string theory, the $g$ loop probability amplitude 
$\mathcal{A}^{(g)}$ is defined as the formal path integral over all two
dimensional compact surfaces with $g$ holes parameterised by $\xi ^{a}$ for $%
a=1,2$ with metric $h_{ab}$ and embedding coordinates $X^{\mu }(\xi ^{1},\xi
^{2})\in \mathbb{R}^{26}$%
\begin{equation}
\mathcal{A}^{(g)}=\frac{1}{V_{\mathrm{Diff}}V_{\mathrm{Weyl}}V_{\mathrm{Trans%
}}}\int \mathcal{D}h\mathcal{D}X\exp (-S_{P}[X,h]),  \label{Ag}
\end{equation}
where 
\begin{equation}
S_{P}[X,h]=\int h^{ab}\partial _{a}X^{\mu }\partial _{b}X^{\mu }\sqrt{\det
(h)}d\xi ^{1}d\xi ^{2}  \label{Spoly}
\end{equation}
is the Polyakov string action \cite{GSW}. $S_{P}$ is invariant under the
group of diffeomorphisms $\xi ^{a}\rightarrow \sigma ^{a}(\xi ^{1},\xi ^{2})$
with formal volume $V_{\mathrm{Diff}}$, the group of Weyl transformations $%
h_{ab}\rightarrow \rho h_{ab}$ with formal volume $V_{\mathrm{Weyl}}$ and
the group of constant translations $X^{\mu }\rightarrow X^{\mu }+a^{\mu }$
with formal volume $V_{\mathtt{Trans}}$. Choose a section (i.e. gauge
fixing) under the diffeomorphism group with complex parameters $z,\bar{z}$
and conformal metric $h_{\mathrm{conf}}\equiv h_{ab}d\xi ^{a}\otimes d\xi
^{b}=e^{\phi }dz\otimes d\bar{z}$ where $z\ $describes a compact Riemann
surface $\mathcal{S}_{g}$ of genus $g$. Formally dividing by $V_{\mathrm{Diff%
}}V_{\mathtt{Weyl}}$ the sum over all metrics reduces to a sum over the
moduli space $\mathcal{M}_{g}$ with Jacobian determinant factors \cite{Kn}.
A standard procedure for doing this is the Faddeev-Popov method where $%
\mathcal{A}^{(g)}$ is re-expressed as a path integral over fermionic ghost $%
b_{zz},c^{z}$and bosonic $X^{\mu }$ configurations \cite{GSW}: 
\begin{equation}
\mathcal{A}^{(g)}=\frac{1}{V_{\mathrm{Trans}}}\int \mathcal{D}b\mathcal{D}c%
\mathcal{D}X\exp (-S_{0}[X]-S_{\mathrm{ghost}}[b,c]),  \label{Aghost}
\end{equation}
with

\begin{eqnarray}
S_{0}[X] &\equiv &S_{P}[X,h_{\mathrm{conf}}]=\int_{\mathcal{S}_{g}}\partial
_{z}X^{\mu }\partial _{\bar{z}}X^{\mu }dzd\bar{z},  \label{S0} \\
S_{\mathrm{ghost}}[b,c] &=&\int_{\mathcal{S}_{g}}(\bar{b}_{\bar{z}\bar{z}%
}\partial _{z}\bar{c}^{\bar{z}}+b_{zz}\partial _{\bar{z}}c^{z})dzd\bar{z}.
\label{Sghost}
\end{eqnarray}
The ghost path integrals reproduce the Jacobian determinants. The bosonic
action $S_{0}$ and the ghost action $S_{\mathrm{ghost}}[b,c]$ are invariant
under conformal transformations of $z$ where $b_{zz}$ transforms as a
conformal 2-form and $c^{z}$ as a 1-vector. The path integral over $X^{\mu }$
results in 
\begin{equation}
\frac{1}{V_{\mathrm{Trans}}}\int \prod_{\mu =1\ldots 26}\mathcal{D}X^{\mu
}\exp (-S_{0}[X])=\left( H_{1}(y,\bar{y})\left| Z_{1}^{(g)}(y)\right|
^{2}\right) ^{26},  \label{Xint}
\end{equation}
where $y\in \mathcal{M}_{g}$ denotes a parameterisation for the complex
moduli for $\mathcal{S}_{g}$. Integration over the constant zero modes
obeying $\partial _{z}X^{\mu }=0$ results in the $H_{1}^{(g)}(y,\bar{y}%
)^{26} $ term on dividing by $V_{\mathrm{Trans}}$. The other term comes from
the path integral over the non-zero modes where $\left|
Z_{1}^{(g)}(y)\right| ^{2}=\det^{-1/2}(-\partial _{z}\partial _{\bar{z}})$
as appropriately defined \cite{Kn}. From the point of view of conformal
field theory, $Z_{1}^{(g)}(y)$ is the genus $g$ partition function for the
one dimensional central charge $C=1$ chiral bosonic CFT (with $\bar{Z}%
_{1}^{(g)}(\bar{y})$ the anti-chiral part). It is also clear that $%
Z_{C}^{(g)}(y)=\left( Z_{1}^{(g)}(y)\right) ^{C}$ where $Z_{C}^{(g)}(y) $ is
the genus $g$ partition function of the $C$ dimensional chiral boson. $%
Z_{C}^{(g)}(y)$ is a function of the moduli $y$ and becomes singular as the
Riemann surface degenerates. Physically, these singularities originate from
the vacuum state of negative energy $-C/24$ where the propagator becomes
singular on pinching down the Riemann surface. For $C=24$, $Z_{24}^{(g)}(y)$
is a meromorphic function of the moduli $y$.

The remaining path integral over the ghost configurations results in the
following factorised form \cite{GSW}, \cite{BK} 
\begin{eqnarray}  \label{Agfact}
\\
\mathcal{A}^{(g)} &=&\int_{\mathcal{M}_{g}}d\mu (y)H_{\mathrm{ghost}%
}^{(g)}(y,\bar{y})\left| Z_{\mathrm{ghost}}^{(g)}(y)\right| ^{2}\left(
H_{1}^{(g)}(y,\bar{y})\left| Z_{1}^{(g)}(y)\right| ^{2}\right) ^{26}, 
\nonumber \\
&=&\int_{\mathcal{M}_{g}}d\mu (y)H^{(g)}(y,\bar{y})\left| F^{(g)}(y)\right|
^{2},  \nonumber
\end{eqnarray}
where 
\begin{eqnarray}
F^{(g)}(y) &=&G^{(g)}(y)Z_{24}^{(g)}(y),  \label{Ffactor} \\
G^{(g)}(y) &=&Z_{\mathrm{ghost}}^{(g)}(y)Z_{2}^{(g)}(y),  \label{Gfactor} \\
H^{(g)}(y,\bar{y}) &=&H_{\mathrm{ghost}}^{(g)}(y,\bar{y})H_{1}^{(g)}(y,\bar{y%
})^{26}.  \label{HFactor}
\end{eqnarray}
Here $\left| Z_{\mathrm{ghost}}^{(g)}(y)\right| ^{2}$ results from
integrating over the ghost non-zero modes where $Z_{\mathrm{ghost}}^{(g)}(y)$
is the genus $g$ partition function for the chiral central charge $-26$
ghost CFT with action $S_{\mathrm{ghost}}[b,c]$. It is a non-singular
function of $y$ since the ghost vacuum energy is $1/12>0$. The zero modes $%
b^{zz}$ are the holomorphic two forms on $\mathcal{S}_{g}$ which are in one
to one correspondence with the moduli \cite{GSW}. Thus the integral measure
for these modes gives us $H_{\mathrm{ghost}}^{(g)}(y,\bar{y})d\mu (y)$ for
some measure $d\mu (y)$ on $\mathcal{M}_{g}.$ Choosing the modular
parameters $y\in \mathcal{T}_{g}$, Teichm\"{u}ller space, we can then take $%
\mu (y)$ to be an appropriate invariant measure under the modular mapping
class group. The integrand $H^{(g)}(y,\bar{y})\left| F^{(g)}(y)\right| ^{2}$
of (\ref{Agfact}) is then a modular invariant.

The meromorphic function $F^{(g)}(y)$ has been discussed by Belavin and
Knizhnik \cite{BK} who argue that $F^{(g)}(y)$ is a meromorphic function
without zeros on $\mathcal{M}_{g}$. They also show that $F^{(g)}(y)$ is
singular only as the Riemann surface degenerates either into two Riemann
surfaces of genus $g_{1},g_{2}$ with $g=g_{1}+g_{2}$ or a single Riemann
surface of genus $g-1,$ as discussed earlier in the genus two case \cite{BK}%
, \cite{Kn}, \cite{Mo}, \cite{KMO}. These singularities are provided by the
vacuum negative energy singularities of $Z_{24}^{(g)}(y)$ so that $%
G^{(g)}(y) $ is a holomorphic function on $\mathcal{M}_{g}$. Furthermore,
the residue of the singular term for $F^{(g)}(y)$ must factorise either into 
$F^{(g_{1})}F^{(g_{2})}$ or $F^{(g-1)}$ depending on the degeneration
considered.

The 1-loop Polyakov string amplitude $\mathcal{A}^{(1)}$ is very well known
where for the standard modular parameter $\tau $ one finds that 
\begin{eqnarray}
Z_{1}^{(1)}(\tau ) &=&\frac{1}{\eta (q)},\ H_{1}^{(1)}(\tau ,\bar{\tau})=%
\frac{1}{\sqrt{Im(\tau )}},  \label{torusdets} \\
Z_{\mathrm{ghost}}^{(1)}(\tau ) &=&\eta (q)^{2},\ H_{\mathrm{ghost}%
}^{(1)}(\tau ,\bar{\tau})=Im(\tau ).  \nonumber
\end{eqnarray}
We then obtain the formal modular invariant amplitude: 
\begin{equation}
\mathcal{A}^{(1)}=\int d\mu (\tau )\frac{1}{Im(\tau )^{12}}\left| \frac{1}{%
\Delta (q)}\right| ^{2},  \label{A1boson}
\end{equation}
with $PSL(2,\mathbb{Z})$ invariant measure $d\mu (\tau )=d^{2}\tau /Im(\tau
)^{2}$. Notice that $G^{(1)}(\tau )=1$ so that $F^{(1)}(\tau
)=Z_{24}^{(1)}(\tau )=1/\Delta (q)$ in this case. Clearly $Z_{24}^{(1)}(\tau
)$ is singular as the torus degenerates to the Riemann sphere for $%
q\rightarrow 0$ as expected.

The 2-loop Polyakov string amplitude $\mathcal{A}^{(2)}$ is much less well
known but has been analysed by Knizhnik \cite{Kn} and others \cite{Mo}, \cite
{KMO}. Then, choosing the period matrix $\Omega $ as moduli, one finds that 
\begin{equation}
H_{1}^{(2)}(\Omega ,\bar{\Omega})=\frac{1}{\sqrt{\det (Im(\Omega ))}},\ H_{%
\mathrm{ghost}}^{(2)}(\Omega ,\bar{\Omega})=\det (Im(\Omega ))^{3},
\label{gentwodets}
\end{equation}
so as to obtain 
\begin{equation}
\mathcal{A}^{(2)}=\int d\mu (\Omega )\frac{1}{\det (Im(\Omega ))^{10}}\left|
F^{(2)}(\Omega )\right| ^{2}.  \label{A2}
\end{equation}
Here the measure $d\mu (\tau )=d^{6}\Omega /\det (Im(\Omega ))^{3}$ is
invariant under the modular group $PSp(4,\mathbb{Z})$ where under a modular
transformation (\ref{modtrans}) $\det (Im(\Omega ))\rightarrow \left| \det
(C\Omega +D)\right| ^{-2}\det (Im(\Omega )).$ Hence from (\ref{A2}), $%
1/F^{(2)}(\Omega )$ transforms as a Siegel 10 form following (\ref{modform})
and is holomorphic because $F^{(2)}(\Omega )$ has no zeros and hence 
\begin{equation}
F^{(2)}(\Omega )={\frac{1}{\Delta _{10}(\Omega )}},  \label{Boson}
\end{equation}
where $\Delta _{10}(\Omega )$ is the Siegel cusp 10 form of (\ref{D10}). The
expected singular factorisation property under $\epsilon \rightarrow 0$
where $F^{(2)}(\Omega )\rightarrow \epsilon
^{-2}F^{(1)}(q_{1})F^{(1)}(q_{2}) $ is then evident from (\ref{D10deg}).
Similarly, we have singular behaviour as $q_{1}\rightarrow 0$ and $%
q_{2}\rightarrow 0.$

In contrast to the genus one case, the holomorphic function $G^{(2)}(\Omega
) $ is non-trivial in the genus two case so that $Z_{24}^{(2)}(\Omega )$ is
not a Siegel modular form. However, we conjecture that they are both
invariant under the translations $T_{1},T_{2}$ and $U$ of (\ref{genset})
whereas 
\begin{eqnarray}
S_{1} &:&Z_{24}^{(2)}(\Omega ){\rightarrow }\frac{{\tau }_{1}^{2}}{{\Omega
_{11}^{12}}}Z_{24}^{(2)}(\Omega ){,}  \label{S1ZB} \\
S_{1} &:&G^{(2)}(\Omega ){\rightarrow }\left( \frac{{\Omega _{11}}}{{\tau }%
_{1}}\right) ^{2}G^{(2)}(\Omega ){,}  \label{S1G}
\end{eqnarray}
and similarly for $S_{2}\,$ so that $F^{(2)}(\Omega )=G^{(2)}(\Omega
)Z_{24}^{(2)}(\Omega )$\thinspace transforms as expected. We will discuss
some evidence to support this conjecture below.

The preceding geometrical analysis can be generalised to 24 dimensional
toroidal compactifications of the bosonic string via a 24 dimensional even
self-dual Niemeier lattice $\Lambda $ where $X^{\mu }\in \mathbb{R}%
^{2}\otimes (\mathbb{R}^{24}/\Lambda ).$ Then the bosonic path integral (\ref
{Xint}) is modified to become: 
\begin{equation}
\frac{1}{V_{Trans}}\int \prod_{\mu =1\ldots 26}\mathcal{D}X^{\mu }\exp
(-S_{0}[X])=\left( H_{1}(y,\bar{y})\left| Z_{1}(y)\right| ^{2}\right)
^{2}\left| Z_{\Lambda }(y)\right| ^{2},
\end{equation}
where the zero modes arise only in the two uncompactified dimensions. $%
Z_{\Lambda }(y)$ is the partition function for the 24 dimensional chiral
lattice theory and is a meromorphic function on $\mathcal{M}_{g}$ with the
same singularities as the uncompactified boson partition function $Z_{24}(y)$%
.

For the genus one case using (\ref{torusdets}) we obtain the standard result 
\begin{equation}
\mathcal{A}_{\Lambda }^{(1)}=\int d\mu (\tau )\left| Z_{\Lambda
}^{(1)}(q)\right| ^{2},  \label{A1lattice}
\end{equation}
where $Z_{\Lambda }(q)$ is modular invariant with the same singularities as $%
1/\Delta (q)$. Hence $Z_{\Lambda }^{(1)}\Delta (q)$ is a holomorphic 12 form
which from CFT is found to be the lattice theta function $\theta _{\Lambda
}(q)$.

For the genus two case using (\ref{gentwodets}) we obtain 
\begin{equation}
\mathcal{A}_{\Lambda }^{(2)}=\int d\mu (\Omega )\det (Im(\Omega ))^{2}\left|
F_{\Lambda }^{(2)}(\Omega )\right| ^{2},  \label{A2lattice}
\end{equation}
with $F_{\Lambda }^{(2)}(\Omega )=G^{(2)}(\Omega )Z_{\Lambda }^{(2)}(\Omega
) $ using the same holomorphic function $G^{(2)}(\Omega )$ as in the
uncompactified bosonic case. From (\ref{A2lattice}), $F_{\Lambda
}^{(2)}(\Omega )$ transforms as a meromorphic Siegel form of weight $2$ and
has the same singularities as $1/\Delta _{10}$. Hence $F_{\Lambda
}^{(2)}(\Omega )\Delta _{10}$ is a holomorphic Siegel 12 form which must be
the theta function $\Theta _{\Lambda }(\Omega )$ from factorisation under
the degeneration $\epsilon \rightarrow 0$. This example is discussed further
below.

\section{Genus Two Bosonic Meromorphic CFT}

\subsection{Bosonic Meromorphic CFT}

We now come to the calculation of the genus two partition function for a
central charge $C$ Bosonic Meromorphic CFT (MCFT) \cite{Go}, \cite{DGM} or
Vertex Operator Algebra \cite{FLM}, \cite{Ka}. We will discuss the examples
of the chiral $C$ dimensional bosonic string, rational MCFTs such as the
even rank $C$ lattice theories and $C=24$ Self-Dual MCFTs which include the
Niemeier lattice theories and orbifoldings such as the Moonshine Module. We
will reproduce the correct singular factorisation and modular properties
described above in string theory from an explicit expression obtained from
sewing together two tori.

Let us briefly discuss some of the features of a bosonic MCFT. Let $\mathcal{%
H}$ be a vector space of states with vacuum state $|0\rangle $ where for
each state $\phi \in \mathcal{H}$ there exists a vertex operator $V(\phi ,z)$
with complex parameter $z\in \mathcal{S}$, the Riemann sphere, with $%
\lim_{z\rightarrow 0}V(\phi ,z)=\phi $ and $V(|0\rangle ,z)=1$. The
operators $\{V(\phi ,z)\}$ satisfy the operator product algebra

\begin{equation}
V(\phi ,z)V(\psi ,w)=V(V(\phi ,z-w)\psi ,w),  \label{VOA}
\end{equation}
with standard ordering $|z|,|z-w|>|w|$ \cite{FLM}, \cite{Go}, \cite{DGM}, 
\cite{Ka}. The automorphism group $G$ of a MCFT is the group of linear
transformations where for $g\in G$ the transformation 
\begin{equation}
gV(\phi ,z)g^{-1}=V(g\phi ,z),  \label{auto}
\end{equation}
preserves (\ref{VOA}). In the case of the Moonshine Module $V^{\natural }$, $%
G\equiv \mathbb{M}$, the Monster finite simple group \cite{FLM}.

The states of $\mathcal{H}$ form a (reducible) representation of the
Virasoro algebra 
\begin{equation}
\lbrack L(m),L(n)]=(m-n)L(m+n)+\frac{C}{12}m(m^{2}-1)\delta (m+n),
\label{Vir}
\end{equation}
where the operator $L(-1)$ generates translations: 
\begin{equation}
\partial _{z}V(\phi ,z)=V(L(-1)\phi ,z).  \label{L-1}
\end{equation}

$\phi $ is said to be quasi-primary if $L(1)\phi =0\,$and primary if
furthermore $L(2)\phi =0$. Each $\phi \in \mathcal{H}$ can be chosen to be
diagonal in $L(0)$ with eigenvalue $h\in \mathbb{Z}_{\ge 0}$, the conformal
weight, where $|0\rangle $ is a primary state of conformal weight $0$. Hence 
$\mathcal{H}$ decomposes into eigenspaces $\mathcal{H}=\oplus _{n\ge 0}%
\mathcal{H}_{n}$. Furthermore, $V(\phi ,z)$ has mode expansion $V(\phi
,z)=\sum_{n\in \mathbb{Z}}\phi (n)z^{-n-h}$ where $\phi (-n):\mathcal{H}%
_{m}\rightarrow \mathcal{H}_{m+n}$ and $\phi =\phi (-h)|0\rangle $ where $%
\phi (-n)$ is said to be of level $n$. As an example, $\mathcal{H\,}$%
contains the quasi-primary Virasoro state $\omega =L(-2)|0\rangle $, of
conformal weight 2, with vertex operator $V(\omega ,z)=\sum_{n\in \mathbb{Z}%
}L(n)z^{-n-2}$ and with modes $L(n)$ obeying (\ref{Vir}).

We also define the Zamolodchikov metric for $\phi ,\psi \in \mathcal{H}$
which is usually heuristically given by the two point function on the
Riemann sphere \cite{BPZ}, \cite{So} 
\begin{equation}
\mathcal{G}_{\phi ,\psi }\equiv \lim_{z_{1}\rightarrow \infty
}\lim_{z_{2}\rightarrow 0}\langle V(\phi ,z_{1})V(\psi ,z_{2})\rangle .
\label{Zam}
\end{equation}
This can be made more rigorous as follows. Under the conformal M\"{o}bius
transformation on the Riemann sphere $z\rightarrow w=-1/z$ a vertex operator
of conformal weight $h_{\phi }$ transforms as \cite{DLM}, \cite{Ka} 
\begin{equation}
V(\phi ,z)\rightarrow w^{2h_{\phi }}V(e^{-L(1)/w}\phi ,w).  \label{mobius}
\end{equation}
We therefore define the Zamolodchikov metric for a MCFT by

\begin{equation}
\mathcal{G}_{\phi ,\psi }|0\rangle \equiv \lim_{w\rightarrow 0}w^{2h_{\phi
}}V(e^{-L(1)/w}\phi ,w)\psi  \label{Zam2}
\end{equation}
One can prove that $\mathcal{G}_{\phi ,\psi }=0$ for $h_{\phi }\neq h_{\psi
} $ and that $\mathcal{G}_{\phi ,\psi }=\mathcal{G}_{\psi ,\phi }$. For
quasi-primary $\phi $ or $\psi $ it is also clear that $\mathcal{G}_{\phi
,\psi }|0\rangle =\phi (h)\psi $. We assume that $\mathcal{G}_{\phi ,\psi }$
is invertible.

Lastly, given two bosonic MCFTs $\{V(\phi ,z)\}$ and $\{V^{\prime }(\phi
^{\prime },z)\}$ with vector spaces $\mathcal{H}$ and $\mathcal{H}^{\prime }$
and central charges $C$ and $C^{\prime }$, one can construct the tensor
product MCFT $\{V(\phi ,z)\otimes V^{\prime }(\phi ^{\prime },z)\}$ with
vector space $\mathcal{H\otimes H}^{\prime }$ and central charge $%
C+C^{\prime }.$ Then the Zamolodchikov metric is 
\begin{equation}
\mathcal{G}_{\phi \otimes \phi ^{\prime },\psi \otimes \psi ^{\prime
}}^{V\otimes V^{\prime }}=\mathcal{G}_{\phi ,\psi }^{V}\mathcal{G}_{\phi
^{\prime },\psi ^{\prime }}^{V^{\prime }}  \label{Zamtensor}
\end{equation}
for $\phi \otimes \phi ^{\prime },\psi \otimes \psi ^{\prime }\in \mathcal{%
H\otimes H}^{\prime }.$

\subsection{Genus One Bosonic MCFT}

The genus two partition function will be constructed by combining
appropriate genus one partition functions. The genus one partition function
is physically defined as the vacuum expectation value on a Riemann torus
with standard parameterisation. For any state $\phi \in \mathcal{H}$, the
torus one point function is defined as the trace \cite{Z}, \cite{So} 
\begin{eqnarray}
Z^{(1)}(\phi ,q) &=&\mathrm{{Tr}_{\mathcal{H}}}(e^{2\pi izh}V(\phi ,e^{2\pi
iz})q^{L(0)-C/24}),  \nonumber \\
&=&\mathrm{{Tr}_{\mathcal{H}}}(\phi (0)q^{L(0)-C/24}),  \label{statetrace}
\end{eqnarray}
where $\phi (0)$ is the level zero mode of $V(\phi ,z)$ (denoted $o(\phi )$
by Zhu and others \cite{Z}, \cite{DLM}). Here  $z\rightarrow e^{2\pi iz}$ is
a conformal transformation from a cylinder (with $z$ and $z+1$ identified)
to the Riemann sphere $\mathcal{S}$. Notice that $Z^{(1)}(\phi ,q)\,$is $z$
independent reflecting the translation symmetry of a torus.

We next define the vertex operator $V[\phi ,z]$ introduced by Zhu \cite{Z} 
\begin{equation}
V[\phi ,z]\equiv e^{2\pi izh}V(\phi ,e^{2\pi iz}-1),  \label{square}
\end{equation}
for $\phi $ of conformal weight $h\ $ (where the origin of the cylinder $z=0$
is mapped into the origin of $\mathcal{S}$). The operators $\{V[\phi ,z]\}$
define a MCFT isomorphic to $\{V(\phi ,z)\}$ with vacuum $|0\rangle $ and
Virasoro state $\tilde{\omega}=\omega -\frac{C}{24}$ $|0\rangle $ with
corresponding modes $\{L[n]\}$ \cite{Z}, \cite{DLM} where in particular 
\begin{eqnarray}
L[0] &=&L(0)+\sum_{n=1}^{\infty }\frac{(-1)^{n+1}}{n(n+1)}L(n),
\label{sqvir} \\
L[-1] &=&L(0)+L(-1).
\end{eqnarray}
$\mathcal{H}$ can thus be alternatively decomposed into eigenspaces of $L[0]$
with $\mathcal{H}=\oplus _{n\ge 0}\mathcal{H}_{[n]}$. As is summarised
below, the torus one point function $Z^{(1)}(\phi ,q)\,$ enjoys simple
modular properties for $\phi $ an eigenstate of $L[0]$.

For the vacuum state of $\mathcal{H}$, (\ref{statetrace}) reduces to the
usual genus one partition function 
\begin{equation}
Z^{(1)}(q)\equiv Z^{(1)}(|0\rangle ,q)=\mathrm{{Tr}_{\mathcal{H}}(}%
q^{L(0)-C/24})\mathrm{.}  \label{Z1}
\end{equation}
For example, in the case of the chiral $C$ dimensional bosonic string, this
is 
\begin{equation}
Z_{C}^{(1)}(q)=\frac{1}{\eta ^{C}(q)},  \label{Boson24}
\end{equation}
whereas for even rank $C$ lattice $\Lambda _{C}$ theories, this is 
\begin{equation}
Z_{\Lambda _{C}}^{(1)}(q)=\frac{\theta _{\Lambda _{C}}(q)}{\eta ^{C}(q)},
\label{Clattice}
\end{equation}
with $\theta _{\Lambda _{C}}(q)=\sum_{\alpha \in \Lambda _{C}}q^{\langle
\alpha ,\alpha \rangle /2}$ which is a modular form for some subgroup of $%
PSL(2,\mathbb{Z})$ of weight $C/2$. For the $C=24$ self-dual MCFTs, the
genus one partition function is modular invariant with 
\begin{equation}
Z_{SD}^{(1)}(q)=\frac{T_{SD}^{(1)}(q)}{\Delta (q)}=J(q)+N_{1},  \label{J}
\end{equation}
where $T_{SD}^{(1)}(q)$ is an elliptic 12 form and $%
J=q^{-1}+0+196884q+O(q^{2})$ is the celebrated modular invariant function
and $N_{1}$ is the number of conformal weight one states. For a MCFT
constructed from a Niemeier lattice $\Lambda $, $T_{SD}^{(1)}=\theta
_{\Lambda }(\tau )$. It is conjectured that there are 71 independent
self-dual MCFTs \cite{Sch}. The Moonshine module $V^{\natural }$ is also
conjectured to be the unique such theory with $N_{1}=0$ \cite{FLM} whereas
of the other 70 conjectured cases, the partition function does not uniquely
determine the CFT \cite{Sch} e.g. $N_{1}=168$ for the lattices $D_{4}^{6}$
and $A_{5}^{4}D_{4}$.

The torus one point function $Z^{(1)}(\phi ,q)$ appears to have received
little attention in the physics literature (with the one possible exception
of \cite{So}) whereas it been extensively discussed by mathematicians such
as Zhu \cite{Z} and Dong, Li and Mason \cite{DLM}, \cite{DM} at least in the
case of rational MCFTs. We note here some of these properties:

\begin{enumerate}
\item[(i)]  For any rational MCFT with $\phi $ of conformal weight $n$ with
respect to $L[0]$, $Z^{(1)}(\phi ,q)$ is a meromorphic elliptic modular form
of weight $n$ \cite{Z}$.$ For the $C$ dimensional boson, $Z_{C}^{(1)}(\phi
,q)$ is an anomalous meromorphic elliptic modular form of weight $n-C$
generalising the behaviour of $\hat{E}_{2}$ in (\ref{SEhat}) \cite{MT}.

\item[(ii)]  Using (\ref{L-1}) we find that for all $\phi $: 
\begin{equation}
Z^{(1)}(L[-1]\phi ,q)=0.  \label{L1Z}
\end{equation}
Geometrically, this follows from the $z$ independence of (\ref{statetrace}).

\item[(iii)]  $Z^{(1)}(\phi ,q)=0$ for all odd $n$ because of the $%
z\rightarrow -z$ reflection symmetry of the torus \cite{Z}, \cite{DLM}.

\item[(iv)]  $Z^{(1)}(\phi ,q)=Z^{(1)}(\phi ^{G},q)$ where $\phi ^{G}$ is
the $G$ invariant part of $\phi $ where $G$ is the automorphism group.
\end{enumerate}

Applying properties (ii)-(iv), it transpires that $Z^{(1)}(\phi ,q)$
vanishes for most states $\phi $ in practice \cite{DM}.

\subsection{Genus Two Bosonic MCFT}

The genus two partition function is physically defined as the vacuum
expectation value of unity on a genus two Riemann surface. This surface can
be constructed by sewing together two tori with modular parameters $%
q_{1},q_{2}$ with a connecting cylinder parameterised by $\epsilon $ as in (%
\ref{pinch}). Then the genus two partition function for a Bosonic MCFT is
explicitly given as:

\begin{equation}
Z^{(2)}(\Omega )=\sum_{n\geq 0}\epsilon ^{n-C/12}\sum_{\phi ,\psi \in 
\mathcal{H}_{[n]}}Z^{(1)}(\phi ,q_{1})\mathcal{G}_{\phi \psi
}^{-1}Z^{(1)}(\psi ,q_{2}),  \label{Z2}
\end{equation}
where the left and right torus one point functions are given in (\ref{Z1}), $%
\mathcal{G}_{\phi \psi }^{-1}$ is the inverse Zamolodchikov metric (\ref
{Zam2}) and the state sum is taken over $\phi ,\psi \in \mathcal{H}_{[n]}$,
the eigenspace of $\mathcal{H}$ for $L[0]$ eigenvalue $n$. We note that
expressions similar to (\ref{Z2}) have long been in existence in the physics
literature e.g. \cite{FS}, \cite{Kn}, \cite{Mo}, \cite{KMO}, \cite{So}, \cite
{P}. However, no completely explicit formula describing the genus two
partition function has formerly appeared with the correct singularity and
modular properties. These are now described.

\textbf{Singularity and Factorising Structure}. $Z^{(2)}(\Omega )$ is a
function of the moduli $q_{1},q_{1},\epsilon $ and is singular and
factorisable under all Riemann surface degenerations. Thus under the
degeneration to two tori with $\epsilon \rightarrow 0$ we have 
\begin{equation}
{Z^{(2)}(\Omega )\rightarrow \epsilon ^{-C/12}Z}^{(1)}(q_{1})Z^{(1)}(q_{2}).
\label{Z2deg}
\end{equation}
Similarly, we have singular factorisation under the other degenerations $%
q_{1}$ (or $q_{2})\rightarrow 0$). It is also clear that the genus two
partition function ${Z^{(2)}(\Omega )}$ for conformal weight $C$ Bosonic
MCFTs has exactly the same singularities as that for the conformal weight $C$
bosonic string. For $C=24$, ${Z^{(2)}(\Omega )\,}$is meromorphic in $%
q_{1},q_{1},\epsilon $.

\textbf{Modular Invariance.} $Z^{(2)}(\Omega )$ does not transform as a
Siegel modular form but has the following modular symmetries. $%
Z^{(2)}(\Omega )$ is obviously invariant under $T_{1},T_{2}$ and $V:\epsilon
\rightarrow -\epsilon $ up to $\frac{24}{C}$th roots of unity. For rational
MCFTs, we find from property (i) above and using (\ref{S1}) that 
\begin{equation}
S_{1}:{Z^{(2)}(\Omega )\rightarrow (-\tau }_{1})^{C/12}{Z^{(2)}(\Omega ).}
\label{S1ZSD}
\end{equation}
and similarly for $S_{2}$. For the $C=24\,$bosonic string, we expect to find
the transformation (\ref{S1ZB}). We will verify this to the first
non-trivial $O(\epsilon ^{2})$ below and note that it has also been verified
to $O(\epsilon ^{6})$ \cite{MT}. For all of the $C=24$ MCFTs, we conjecture
that ${Z^{(2)}(\Omega )}$ is also invariant under the $U$ modular
transformation. This should follow from re-expressing (\ref{Z2}) as a sewn
torus two point function together with the associativity property of the
operator product expansion (\ref{VOA}) \cite{MT}.

Multiplying ${Z_{24}^{(2)}(\Omega )\,}$by the holomorphic correction factor $%
{G}^{(2)}{(\Omega )}$ of (\ref{Gfactor}), we obtain a Siegel Modular form ${F%
}^{(2)}{(\Omega )=G}^{(2)}{(\Omega )Z_{24}^{(2)}(\Omega )}$ where from (\ref
{S1ZB}) and (\ref{S1G}) we see that ${F}^{(2)}{(\Omega )}$ transforms as a
meromorphic Siegel modular form of weight $-10$ so that ${F}^{(2)}{(\Omega
)=\Delta }_{10}^{-1}(\Omega )$ as discussed earlier.

Consider next a $C=24$ self-dual MCFT with genus two self-dual partition
function ${Z_{SD}^{(2)}(\Omega )}$ which has the same singularities as ${%
Z_{24}^{(2)}(\Omega ).}$ Then, ${F}_{SD}^{(2)}{(\Omega )=G}^{(2)}{(\Omega
)Z_{SD}^{(2)}(\Omega )}$ transforms as a meromorphic Siegel modular form of
weight $2$ with the same singularities as ${\Delta }_{10}^{-1}(\Omega )$
where we use the same holomorphic correction factor. Equivalently, if we
write 
\begin{equation}
{Z_{SD}^{(2)}(\Omega )=}T_{SD}^{(2)}({\Omega }){Z_{24}^{(2)}(\Omega ),}
\label{T2}
\end{equation}
then $T_{SD}^{(2)}({\Omega })$ is a holomorphic Siegel modular form of
weight $12$. From (\ref{Z2}) this factorises under the degeneration $%
\epsilon \rightarrow 0$ with 
\begin{equation}
T_{SD}^{(2)}(\Omega )\rightarrow T_{SD}^{(1)}(q_{1})T_{SD}^{(1)}(q_{2}).
\label{T2deg}
\end{equation}
with $T_{SD}^{(1)}(q)$ as in (\ref{J}). However, every Siegel 12 form is
uniquely determined by this degeneration. In particular, for the Niemeier
lattice theories with $T_{SD}^{(1)}(q)=\theta _{\Lambda }(q)$ we therefore
find that $T_{SD}^{(2)}(\Omega )=\Theta _{\Lambda }(\Omega )$, the Siegel
lattice theta function (\ref{ThetaL}) and hence 
\begin{equation}
{Z_{SD}^{(2)}(\Omega )=\Theta _{\Lambda }(\Omega )Z_{24}^{(2)}(\Omega ),}
\end{equation}
in these cases.

In general, we may express $T_{SD}^{(2)}(\Omega )$ in terms of the basis of
Siegel 12 forms $\psi _{4}^{3}$, $\psi _{6}^{2}$ and $F_{12}$ to find the
following explicit formula 
\begin{equation}
T_{SD}^{(2)}(\Omega )=c_{1}\psi _{4}^{3}+c_{2}\psi
_{6}^{2}+(1-c_{1}-c_{2})F_{12},  \label{T2formula}
\end{equation}
where $c_{1}=(1927+6k-k^{2})/1152$, $c_{2}=(1457-78k+k^{2})/6336$ and $%
k=N_{1}/24$ is the dual Coxeter number in the Niemeier lattice cases whereas 
$k=1$ for the Leech lattice theory $V^{\Lambda }$ and $k=0$ for the
Moonshine module $V^{\natural }$ with Fourier expansion 
\begin{equation}
T^{\natural (2)}=(1-24q)(1-24s)+48qsu-24qsu^{2}+O(q^{2},s^{2}).
\label{TMoon}
\end{equation}
Note that the general formula (\ref{FkFourier}) is again satisfied.

A surprising consequence of the general formula (\ref{T2formula}) is that
the genus two partition function does not always distinguish between all the
conjectured 71 $C=24$ self-dual MCFTs since it is completely determined by
its singular factorisation under $\epsilon \rightarrow 0$. Thus if two $C=24$
self-dual MCFTs share the same genus one partition function, they share the
same partition function at genus two also.

\textbf{Tensor products of MCFTs.} If $\{V(\phi ,z)\}$ and $\{V^{\prime
}(\phi ^{\prime },z)\}$ are two MCFTs with genus two partition functions ${%
Z_{V}^{(2)}(\Omega )}$ and ${Z_{V^{\prime }}^{(2)}(\Omega )}$ respectively,
then it is clear using (\ref{Zamtensor}) that the tensor product $\{V(\phi
,z)\otimes V^{\prime }(\phi ^{\prime },z)\}$ MCFT has genus two partition
function 
\begin{equation}
{Z_{V\otimes V^{\prime }}^{(2)}(\Omega )=Z_{V}^{(2)}(\Omega )Z_{V^{\prime
}}^{(2)}(\Omega ),}  \label{Ztensor}
\end{equation}
Thus the dimension $C$ bosonic string MCFT is the tensor product of $C$
copies of the $C=1$ bosonic string MCFTs so that 
\begin{equation}
{Z_{C}^{(2)}(\Omega )=}\left( {Z_{1}^{(2)}(\Omega )}\right) ^{C},
\label{ZBosonC}
\end{equation}
as expected from the path integral analysis earlier.

It is not difficult to see that much of the previous analysis carries
directly over to other rational MCFTs such as for an even rank $C$ lattice $%
\Lambda _{C}$ MCFT. Then (\ref{T2}) generalises so that 
\begin{equation}
{Z_{\Lambda _{C}}^{(2)}(\Omega )=}\Theta _{\Lambda _{C}}(\Omega ){%
Z_{C}^{(2)}(\Omega ),}  \label{TC}
\end{equation}
where $\Theta _{\Lambda _{C}}(\Omega )$ is the Siegel lattice $\Lambda _{C}$
theta function which is a modular form of weight $C/2$ for some appropriate
subgroup of $PSp(4,\mathbb{Z})$ in general. As noted earlier, for $C=0\text{
mod }8$, $\Theta _{\Lambda _{C}}(\Omega )$ is a modular form of weight $C/2$
for the full $PSp(4,\mathbb{Z})$ group.

\subsection{Some Examples}

We will discuss a number of examples of genus two partition functions (\ref
{Z2}) expanded to $O(\epsilon ^{2})$. Many of these examples have also been
computed and confirmed to have the correct modular properties to $O(\epsilon
^{6})$ \cite{MT}. For the $C$ dimensional bosonic string, even rank lattice
theories and $C=24$ self-dual MCFTs, we can use properties (ii)-(iv) above
to find that the only non-zero trace that arises at $L[0]$ conformal weight
two is that for the Virasoro state $\tilde{\omega}=\omega -C/24$ where \cite
{DLM}, \cite{DM}: 
\begin{eqnarray}
Z^{(1)}(\tilde{\omega},q) &=&\mathrm{{Tr}_{\mathcal{H}}}(\tilde{\omega}%
(0)q^{L(0)-C/24}),  \label{Zw} \\
&=&\mathrm{{Tr}_{\mathcal{H}}}((L(0)-\frac{C}{24})q^{L(0)-C/24}),  \nonumber
\\
&=&q\frac{d}{dq}Z^{(1)}(q).  \nonumber
\end{eqnarray}
For the $C$ dimensional bosonic string we obtain an anomalous $2-C$ form 
\begin{equation}
Z_{C}^{(1)}(\tilde{\omega},q)=\frac{C}{2}\frac{\hat{E}_{2}(q)}{\eta ^{C}(q)},
\label{Zwboson}
\end{equation}
whereas for the lattice $\Lambda _{C}$ theories and the $C=24$ self-dual
MCFTs we obtain 2-forms 
\begin{eqnarray}
Z_{\Lambda _{C}}^{(1)}(\tilde{\omega},q) &=&\frac{D_{q}\theta _{\Lambda
_{C}}(q)}{\eta ^{C}(q)},  \label{Zwlattice} \\
Z_{SD}^{(1)}(\tilde{\omega},q) &=&\frac{D_{q}T^{(1)}(q)}{\Delta (q)}.
\label{ZwSD}
\end{eqnarray}
Here we are using the covariant derivative $D_{q}$ of (\ref{Df}) and the
fact that $q\frac{d}{dq}\eta =-\frac{1}{2}\hat{E}_{2}\eta $ (since $%
D_{q}\Delta =0$ because there are no elliptic cusp forms of weight 14)$.$

The state $\tilde{\omega}$ is quasi-primary so we have $\mathcal{G}_{\tilde{%
\omega},\tilde{\omega}}|0\rangle =L[2]\tilde{\omega}=\frac{C}{2}|0\rangle $
from (\ref{Vir}) $\,$i.e. $\mathcal{G}_{\tilde{\omega},\tilde{\omega}}^{-1}=%
\frac{2}{C}$. Hence from (\ref{Z2}) we obtain the genus two partition
function to $O(\epsilon ^{2})$ to be 
\begin{equation}
Z^{(2)}(\Omega )=\frac{1}{\epsilon ^{C/12}}\left(
Z^{(1)}(q_{1})Z^{(1)}(q_{2})+\frac{2}{C}Z^{(1)}(\tilde{\omega},q_{1})Z^{(1)}(%
\tilde{\omega},q_{2})\epsilon ^{2}+O(\epsilon ^{4})\right) .  \label{Z2eps}
\end{equation}
Hence we find for the $C$ dimensional boson that 
\begin{equation}
Z_{C}^{(2)}(\Omega )=\frac{1}{\epsilon ^{C/12}}\frac{1}{\eta ^{C}(q_{1})}%
\frac{1}{\eta ^{C}(q_{2})}\left( 1+\frac{C}{2}\hat{E}_{2}(q_{1})\hat{E}%
_{2}(q_{2})\epsilon ^{2}+O(\epsilon ^{4})\right) .  \label{Z2BosonC}
\end{equation}
It is easy to confirm the modular property (\ref{S1ZB}) and the product
formula (\ref{Z2BosonC}) to the given order in $\epsilon ^{2}$ by using (\ref
{SEhat}) and (\ref{om11}). These have also been confirmed for expansions up
to $O(\epsilon ^{6})$ \cite{MT}.

The explicit calculation of the universal holomorphic function $%
G^{(2)}(\Omega )\,$of (\ref{Gfactor}) from a sewing formalism requires a
formulation of (\ref{Z2}) for the fermionic central charge $-26$ ghost CFT
with negative Zamolodchikov metric elements. This is currently under
investigation \cite{MT} but we conjecture that the genus two ghost partition
function is 
\begin{equation}
Z_{\mathrm{ghost}}^{(2)}(\Omega )=\epsilon ^{1/6}\eta ^{2}(q_{1})\eta
^{2}(q_{2})\left( 1-3\hat{E}_{2}(q_{1})\hat{E}_{2}(q_{2})\epsilon
^{2}+O(\epsilon ^{4})\right) .  \label{Z2ghost}
\end{equation}
We then obtain from (\ref{Ffactor}) and (\ref{Gfactor}) that 
\begin{eqnarray}
\Delta _{10}(\Omega ) &=&\epsilon ^{2}\Delta (q_{1})\Delta (q_{2})\left( 1-10%
\hat{E}_{2}(q_{1})\hat{E}_{2}(q_{2})\epsilon ^{2}+O(\epsilon ^{4})\right) ,
\label{D102} \\
G^{(2)}(\Omega ) &=&1-2\hat{E}_{2}(q_{1})\hat{E}_{2}(q_{2})\epsilon
^{2}+O(\epsilon ^{4}).  \label{G2}
\end{eqnarray}
It is easy to confirm the novel $\epsilon $ expansion for the Siegel 10 form 
$\Delta _{10}(\Omega )$ and that $G^{(2)}(\Omega )\,$ has the modular
property (\ref{S1G}) to the given order in $\epsilon ^{2}$ by using (\ref
{SEhat}) and (\ref{om11}). Using the formulas (\ref{qsueps}) for the Fourier
parameters $q,s,u$ and that $\hat{E}_{2}(q)=-\frac{1}{12}+2q+O(q^{2})$ we
may expand (\ref{D102}) to $O(q_{1}^{2},q_{2}^{2},\epsilon ^{4})$ to also
confirm the Fourier expansion (\ref{D10fourier}).

The even rank lattice and $C=24$ self-dual MCFTs have genus two partition
functions 
\begin{eqnarray}
Z_{\Lambda _{C}}^{(2)}(\Omega ) &=&\frac{1}{\epsilon ^{C/12}}\frac{\theta
_{\Lambda _{C}}(q_{1})}{\eta ^{C}(q_{1})}\frac{\theta _{\Lambda _{C}}(q_{2})%
}{\eta ^{C}(q_{2})}(1+  \label{Z2latticeeps} \\
&&\frac{2}{C}\frac{D_{q}\theta _{\Lambda _{C}}(q_{1})}{\theta _{\Lambda
_{C}}(q_{1})}\frac{D_{q}\theta _{\Lambda _{C}}(q_{2})}{\theta _{\Lambda
_{C}}(q_{2})}\epsilon ^{2}+O(\epsilon ^{4})),  \nonumber \\
Z_{SD}^{(2)}(\Omega ) &=&\frac{1}{\epsilon ^{2}}\frac{T^{(1)}(q_{1})}{\Delta
(q_{1})}\frac{T^{(1)}(q_{2})}{\Delta (q_{2})}(1+  \label{Z2SDeps} \\
&&\frac{1}{12}\frac{D_{q}T^{(1)}(q_{1})}{T^{(1)}(q_{1})}\frac{%
D_{q}T^{(1)}(q_{2})}{T^{(1)}(q_{2})}\epsilon ^{2}+O(\epsilon ^{4})). 
\nonumber
\end{eqnarray}
\newline

The Siegel forms $T^{(2)}(\Omega )$ of (\ref{T2}) and $\Theta _{\Lambda
_{C}}(\Omega )$ of (\ref{TC}) for $C=0\text{ mod }8$, lead us to the
following novel $\epsilon $ expansion for a Siegel $k=C/2$ form $%
F_{k}(\Omega )$ which factorises $F_{k}(\Omega )\rightarrow
f_{k}(q_{1})f_{k}(q_{2})$ as $\epsilon \rightarrow 0$ (where $f_{k}(q)$ is
an elliptic $k$ form) as follows: 
\begin{eqnarray}
F_{k}(\Omega )&=&f_{k}(q_{1})f_{k}(q_{2})(1+  \label{Fkeps} \\
&&(\frac{1}{k}\frac{D_{q}f_{k}(q_{1})}{f_{k}(q_{1})}\frac{D_{q}f_{k}(q_{2})}{%
f_{k}(q_{2})}-k\hat{E}_{2}(q_{1})\hat{E}_{2}(q_{2}))\epsilon ^{2}+O(\epsilon
^{4})).  \nonumber
\end{eqnarray}
For $f_{k}(q)=1+aq+O(q^{2})$ and using (\ref{qsueps}) we may expand (\ref
{Fkeps}) to $O(q_{1}^{2},q_{2}^{2},\epsilon ^{4})$ to confirm the expression
(\ref{FkFourier}) for $C=0\text{ mod }8$.

\section{Future Directions}

We have described the genus two partition function for bosonic Meromorphic
CFTs by sewing together appropriate torus one point functions. A number of
important issues remain outstanding:

\begin{itemize}
\item  Can the modularity property (\ref{S1ZB}) for $Z_{24}^{(2)}(\Omega )$
be rigorously established? Likewise, $U$ invariance needs to be shown.

\item  The universal holomorphic correction $G^{(2)}(\Omega )$ which
contains the ghost contributions needs to be understood. Can we calculate
and motivate the need for $Z_{\mathrm{ghost}}^{(2)}(\Omega )\,$ within a
MCFT/VOA formulation?

\item  What, if any, is the number theoretic meaning of $Z_{1}^{(2)}(\Omega )
$ and can a closed form be found for it?
\end{itemize}

Finally, amongst the many future possible developments we mention genus two
orbifold constructions. For a pair of elements $g_{1},g_{2}\in G$, the
automorphism group, we can define \newline

\begin{equation}
Z^{(2)}(g_{1},g_{2},\Omega )=\sum_{n\geq 0}\epsilon ^{n-C/12}\sum_{\phi
,\psi \in \mathcal{H}_{[n]}}Z^{(1)}(g_{1},\phi ,q_{1})\mathcal{G}_{\phi \psi
}^{-1}Z^{(1)}(g_{2},\psi ,q_{2}),  \label{Z2orb}
\end{equation}
where 
\begin{equation}
Z^{(1)}(g,\phi ,q)=\mathrm{{Tr}_{\mathcal{H}}}(g\phi (0)q^{L(0)-C/24}),
\end{equation}
as discussed by Dong, Li and Mason \cite{DLM}, \cite{DM}. The modular
properties for (\ref{Z2orb}) would be of particular interest where we expect
that twisted sector traces will appear. In the case of the Moonshine Module,
we should in principle be able to obtain every Generalised Moonshine trace
function \cite{N}, \cite{T}, \cite{DLM} by choosing some commuting elements $%
g_{1},g_{2}$ on a genus two Riemann surface, performing an appropriate
modular transformation and then pinching down to a genus one surface with
twisted boundary conditions. In this way, we could calculate Generalised
Moonshine functions directly from the Moonshine Module and hopefully,
understand the genus zero property in a new way.

\end{document}